\documentclass[12pt]{amsart}
\usepackage{amssymb}
\allowdisplaybreaks[1]

\theoremstyle{definition}

\numberwithin{equation}{section}

\newcommand{\R}{\mathbb R}

\renewcommand{\a}{\alpha}

\renewcommand{\v}{\varphi}

\newcommand{\s}{\sigma}

\newcommand{\e}{\xi}

\renewcommand{\d}{\delta}

\renewcommand{\l}{\lambda}
\font\got=eufm10
\renewcommand{\u}[1]{\hbox{\got #1}}

\newcommand{\p}{\partial}

\renewcommand{\k}{\prime}

\begin{document}

\title[Tangential Star Products]
{Tangential Star Products}

\author[\small Angela Gammella]{A. Gammella}
\address{
UNIVERSITE LIBRE DE BRUXELLES\\
CAMPUS DE LA PLAINE\\
CP 218, BOULEVARD DU TRIOMPHE\\
1050 BRUXELLES, BELGIQUE.}
\email{agammell@ulb.ac.be}

%---------------------------------------------------------------------

\

\maketitle

\

\

\begin{abstract}
  We etablish a necessary and sufficient condition under which there exists a tangential
and well graded 
 star product, differential or not, on the dual $\u{g}^\ast$ of a nilpotent Lie algebra $\u{g}$. We also
give  enlightening examples with explicit computations.
\end{abstract}

\

\

\

\

\centerline{\bf Introduction}

\

The theory of deformations and especially the notion of star products have been developed by Flato,
Lichnerowicz and  their collaborators in [3] with the aim of quantizing a classical system represented by
a symplectic or a Poisson  manifold $M$.
 \noindent
A star product on $M$ is an associative, non commutative product on
$C^{\infty}(M)$ depending formally on a parameter 
$\nu$ (in physical applications $\nu$ is $\frac{i \hbar}{2}$ where $\hbar $ denotes Planck's
constant).  The product should have the form
 $$\displaystyle f\ast g = \sum_{n\geq 0} C_n(f,g)\nu ^n $$
\noindent where $f,\,g$ are in $C^{\infty}(M)$, $C_0(f,g)=fg$, $C_1(f,g)=\{f,g\}$ and $C_n(f,g)$ are bilinear operators 
on $C^{\infty}(M)$ with values in $C^{\infty}(M)$.\par
\noindent
The main development of this theory went through the proof of the existence of differential 
star products, that
is star products whose cochains $C_n$ are differential operators. In the case of symplectic manifolds, the
question has been completely solved, using different approaches [7, 8]. In a recent work, Kontsevich has
given a remarkable proof of the existence of differential star products on an arbitrary smooth Poisson
manifold [10].

\par
 Since every Poisson manifold is ``foliated" by symplectic submanifolds, it is quite natural to study 
 star products with
 nice restrictions to the symplectic leaves. Such  star products are called tangential.
 For regular Poisson manifolds, a tangential version of Vey's work is introduced in [11] and a proof of the existence of
 tangential and differential star products can be found in [12]. Unfortunately, for general Poisson
manifolds,
 such tangential and differential star products do not always exist.

\noindent
 Indeed, let $\u{g}$ be a Lie algebra and let us consider the dual space $\u{g}^\ast $, endowed with its linear
 Poisson structure. It is well known that in this case,
 the symplectic leaves are nothing else but the coadjoint orbits in $\u{g}^{\ast}$. 
 It turns out that a tangential and differential deformation on $\u{g}^\ast$ is possible 
 only if $\u{g}$ satisfies a very strong algebraic condition [5]. No semi-simple Lie algebras satisfy this condition. 
 Moreover, it has been shown that the standard deformation on $\u{g}^{\ast}$, {\sl i.e.} the Gutt 
 star product, is very
 rarely tangential [2].

\par
However, for non-differential star products, the situation is far better. Cahen and Gutt have
constructed in [6] an algebraic tangential star product on the set of regular orbits of any semi-simple
Lie algebra.  Furthermore, there is in [1] a construction of a deformation on $\u{g}^\ast$ in the case
where the symmetric algebra
$S(\u{g})$ is a free $I(\u{g})$-module ($I(\u{g})$ denoting the algebra of invariant polynomials on $\u{g}^\ast$). This
deformation generalizes the one which is given in the semi-simple case in [6], but is less explicit.

\par
The purpose of the present paper is to give a simple condition, Theorem 10, for the existence of
tangential star products  on $\u{g}^\ast$ in the nilpotent case. 
First, we introduce the needed notions and compute the cohomology related to deformations of the associative and
graded algebra $S(\u{g})$. Then, we prove that the construction of a good operator $C_2$ is enough to ensure the 
existence of a tangential, ``well"  graded, differential or not, star product on 
$\u{g}^{\ast}$. We devote the last part to explicit illustrations. In particular, 
we apply our result to $\u{g}_{54},$ the simplest example of a nilpotent Lie algebra $\u{g}$ for
which there is no tangential and differential deformation of $S(\u{g})$ [1, 5]. 

\medskip{\bf Notation:}
Throughout this letter, $\u{g}$ denotes a nilpotent Lie algebra and $\u{g}^\ast$ the dual space of
$\u{g}$.  The symmetric algebra $S(\u{g})$ over $\u{g}$ is naturally identified with the algebra of
real-valued polynomials on the  dual $\u{g}^\ast$. Obviously, $S(\u{g})=\oplus S^k(\u{g})$ where
$S^k(\u{g})$ is the space of homogeneous polynomials  of degree $k$. We denote by $I(\u{g})$ (or $I$) the
algebra of invariant polynomials on $\u{g}^\ast$.
\vskip1cm
\section{Differential, algebraic and tangential operators}\label{sec1}
\smallskip We first recall some essential facts about nilpotent Lie algebras. See [1, 4, 15] for
more details.\par
\noindent
Suppose that $\u{g}$ is an $m$-dimensional nilpotent Lie algebra. Denote by $(X_i)$ a Jordan-H\"older
basis of $\u{g}$ (that is $[X_i,X_j]\equiv 0\,\hbox{mod}\,(X_1,\ldots,X_{j-1})$  if $i\geq j$).
Let $(x_i)$ be the system of coordinates of $\u{g}^\ast$ associated to this basis.
 Let $G$ be the simply connected group with Lie algebra $\u{g}$. Let also $2d$ be the maximal dimension
of coadjoint orbits in $\u{g}^\ast$. There exist: \par
(i) a Zariski open subset $V$ of $\u{g}^\ast$, invariant by the action of the adjoint group of $G$, dense in $\u{g}^\ast$,
containing only orbits of maximal dimension;\par
(ii) $2d$ rational functions $(p_1,...,p_d,q_1,...q_d)$ in the variables $(x_i)$ which are regular on
$V$;\par (iii) $m-2d$ polynomial functions ${\l}_1,...,\l_{m-2d}$ in the variables $(x_i)$;\par
(iv) a Zariski open subset $U$ of $\R^{m-2d}$.\par\noindent
These elements are such that there exists
a diffeomorphism $\v$ between $V$ and $U\times \R^{2d}$ defined
by $\v(\e)=(\l(\e),p(\e),q(\e))$ if we note $\l=({\l}_1,...,\l_{m-2d})$, $p=(p_1,...,p_d)$ and 
$q=(q_1,...,q_d)$, such that each orbit contained 
in $V$ admits a global Darboux chart defined by the variables $p_i,q_j$ and that each invariant rational function
on $\u{g}^\ast$ may be written in a unique way as a rational function in the variables $({\l}_k)$. 
The orbits contained in $V$ are usually called generic orbits and each polynomial $\l_k$ is said to be a generic invariant.
Moreover, every $X$ in $\u{g}$, as a function on $\u{g}^\ast$, restricted to $V$, can be written as
$$X=\displaystyle \sum_{1\leq j\leq d} a_j(q,\l) p_j + a_0(q,\l)$$
where the coefficients $a_j$ are polynomial in $q$ and rational in $\l$. \par
 Let us denote by $\R(\l)[p,q]$ the algebra of polynomial functions in $p,q$ with coefficients in the space $\R(\l)$ of 
rational functions in $\l$. Thus, every $X$ in $\u{g}$ can be identified with an element of $\R(\l)[p,q]$. 
Let us now consider $S(\u{g})_I$, the localized algebra of rational functions with non zero invariant denominators.
We see that the quotient field of $I$ is exactly $\R(\l)$, thus $S(\u{g})_I$ coincides with 
the space $\R(\l)[x_1,...,x_m]$ of polynomials on $\u{g}^\ast$ with rational coefficients in $\l$ and also
with
$\R(\l)[p,q]$. Furthermore, for each element $X$ of $\u{g}$ considered as a function on $\u{g}^\ast$, the
derivative $\p_{X}$ with respect  to $X$ can be written in the form
$$\p_X=\displaystyle \sum_{1\leq i\leq d} a_i\p_{p_i} +\sum_{1\leq j\leq d} b_j \p_{q_j} + \sum_{1\leq
k\leq m-2d} c_k \p_{{\l}_k}$$
with $a_i,\,b_j$ and $c_k$ in $S(\u{g})_I\simeq\R(\l)[x_1,...,x_m]\simeq\R(\l)[p,q]$.\medskip\par

Let us then fix the system of linear coordinates $(x_i)$ on $\u{g}^\ast$ and the local system of coordinates $(p,q,\l)$
as above. From now on, the localized algebra $S(\u{g})_I$ will
be identified with $\R(\l)[p,q]$.\par

\medskip\noindent
{\bf Definition 1.}\par\noindent
{\sl A multilinear map $F: S(\u{g})\times\ldots\times S(\u{g})\rightarrow S(\u{g})$ is said to be
differential if it is given by differential operators ({\sl i.e.} of finite order) on each argument.
Otherwise,  it is called algebraic or non differential.}\par\smallskip 
Now, let $D$ be a differential operator on $S(\u{g})$. Then, we can write 
$$D(u_1,...,u_s)=\sum D_{{\a}_1...\a_s} \p_{{\a}_1}(u_1) \ldots \p_{{\a}_s}(u_s)$$ where the multi-indexes
${\a}_i$ are relative to the variables $(x_i)$ and where the coefficients $D_{{\a}_1...\a_s}$ belong to
$S(\u{g})$. The same operator can be  written as a differential operator, say $\widetilde{D}$, in the
variables $(p,q,\l)$ just by performing, for the  coefficients and for the operators ${\p}_{\a_i}$, the
change of variables from $(x_i)$ to $(p,q,\l)$. Such a differential  operator $D$ (or ${\widetilde{D}}$) can
naturally be extended to the localized algebra $S(\u{g})_I$.
  Now, let $A$ be an algebraic operator on $S(\u{g})$. $A$ can be decomposed into an infinite sum 
$\displaystyle\sum_N A_N$ of differential operators of the form 
$$A_N(u_1,...,u_s)=\displaystyle\sum_{|\a_1|+\ldots+|\a_s|=N} A_{{\a}_1...\a_s} \p_{{\a}_1}(u_1) \ldots
\p_{{\a}_s}(u_s)$$ where $\a_i=(\a{}_i^1,\ldots,\a{}_i^m)$ are multi-indexes relative to the variables
$(x_i)$ and 
$|\a_i|=\a{}_i^1+\ldots +\a{}_i^m.$  Let ${\widetilde{A}}_{(t)}$ be the operator defined by 
${\widetilde{A}}_{(t)}=\sum t^N \widetilde{A_N}$. Clearly,
${\widetilde{A}}_{(1)}$ coincides with $A$ on $S(\u{g})$ and ${\widetilde{A}}_{(t)}$ sends $S(\u{g})_I$ into the (formal) algebra 
$\R(\l)[p,q][[t]]$ of formal series in $t$ with coefficients in $\R(\l)[p,q]$.
\medskip
In the following, we shall use the algebraic notion of tangential operators given in [5].

\medskip\noindent
{\bf Definition 2.}\par\noindent
{\sl A multilinear map $F:S(\u{g})\times ...\times S(\u{g})\rightarrow
S(\u{g})$ is 
said to be
tangential if $F$ vanishes on constants and if for each $\Delta$ in $I=I(\u{g})$, for every 
$u_1,...,u_s$ in $S(\u{g})$
and for all $1\leq l\leq s$, $$\Delta F(u_1,...,u_s)=
F(u_1,...,\Delta\,u_l,...,u_s).$$}
\par\smallskip\noindent Such a tangential operator can be uniquely
extended to the localized algebra $S(\u{g})_I$ of rational functions with non zero invariant denominators
by
$$\hat{F}(\frac{u_1}{Q_1},\ldots,\frac{u_s}{Q_s})=\frac{1}{Q_1...Q_s} F(u_1,...,u_s)$$
where the $u_i$ are in $S(\u{g})$ and the $Q_i$ are elements of $I$.\par
Now, it is possible to characterize tangential operators thanks to the variables $(p,q,\l)$. 
Indeed, if $F$ is
a tangential map on $S(\u{g})$, then its extension $\hat{F}$ to $S(\u{g})_I$ satisfies
$$\hat{F}(v_1,\ldots,\l_k v_l,\ldots,v_s)=\l_k \hat{F}(v_1,...,v_s)$$
for all generic invariant $\l_k$ and for all $v_i$ in $S(\u{g})$. It follows that $\hat{F}$ is of the form
$$\hat{F}(v_1,...,v_s)=\sum F_{{\widetilde{\a}_1} \ldots
{\widetilde{\a}_s}}(p,q,\l)\p_{{\widetilde{\a}_1}}(v_1)\ldots\p_{{\widetilde{\a}_s}}(v_s)$$ here the
${\widetilde{\a}_i}$ are multi-indexes relative to the variables $(p,q,\l)$, the coefficients 
$F_{{\widetilde{\a}_1} \ldots {\widetilde{\a}_s}}$ belong to $S(\u{g})_I\simeq \R(\l)[p,q]$ and the
$\p_{{\widetilde{\a}_i}}$ do not  include derivatives with respect to the variables $(\l_k)$.
\par\noindent
Conversely, suppose that $C=\sum C_N$ is an algebraic operator on $S(\u{g})$ such that 
${\widetilde{C}}_{(1)}$
can be expressed without derivatives with respect to $(\l_k)$, then $C$ is tangential.
\vskip0,5cm\noindent
{\bf Remark 3.}\par\noindent
{\sl  Frequently, a tangential operator $F$ on $S(\u{g})$,  
in the sense of Definition 2, is not only tangential on the set $V$ of generic orbits 
but also on the set $\Omega$ of all orbits of maximal dimension (see the example of 
$\u{g}_{54}$ in Section 4
for instance).
More precisely, each tangential mapping $F$ is in fact tangential on any open subset $O$ of 
$\u{g}^\ast$, such that the only polynomials on $\u{g}^\ast$, whose restriction to an orbit contained in $O$ is
zero, belong to $I$. This set $O$ contains $V$ by construction and often coincides with $\Omega$.}    

\vskip1cm

\section{Well graded cohomology}\label{sec2}

\smallskip\noindent
{\bf Definition 4.}\par\noindent
{\sl
A $s$-linear map $C:S(\u{g})\times\ldots\times S(\u{g}) \rightarrow S(\u{g}) $ is said to be 
homogeneous of degree $-n$ 
if for all $u_1,...,u_s$ in $S^{d_1}(\u{g}),\ldots,S^{d_s}(\u{g})$  respectively, $C(u_1,...,u_s)$ is in
 $S^{d_1+\ldots+d_s-n}(\u{g})$.}\par \smallskip
Now, let us recall that each differential operator on $S(\u{g})$ can naturally be extended to $S(\u{g})_I$. If $u$ belongs 
to $S(\u{g})_I\simeq \R(\l)[p,q]$, we will denote by $||u||$ the degree of $u$ as a polynomial in $p$.
Then, the following definition makes sense.\par
\medskip\noindent
{\bf Definition 5.}\par\noindent
{\sl A $s$-differential operator $D:S(\u{g})\times\ldots\times S(\u{g})\rightarrow S(\u{g})$ is said to
be 
correct of degree $-n$
if for all $u_i$ in $S(\u{g})_I$ such that $||u_i||=d_i$, $(1\leq i\leq s)$, 
$$||D(u_1,...,u_s)||\leq (d_1+\ldots+d_s)-n.$$
An algebraic operator $C=\sum C_N$ on $S(\u{g})$ is said to be correct of degree $-n$ if all the 
differential operators $C_N$ are correct of degree $-n$.} 

\vskip0,5cm
Let ${\mathcal A}$ be a commutative and associative algebra and $M$ be a ${\mathcal A}$-bimodule. We can introduce the graded 
${\mathcal A}$-module of Hochschild cochains $C^{\ast}({\mathcal A},M)$ that is the ${\mathcal A}$-module of multilinear 
maps with values in $M$. \par
\medskip
\noindent
{\bf Definition 6.}\par\noindent
{\sl The coboundary of a $s$-cochain $C:{\mathcal A}\times\ldots\times{\mathcal A} \longrightarrow M$ is
the
$(s+1)$-cochain $\d C$  defined by
\begin{align*}
\d C(u_1,...,u_{s+1}) & =u_1 C(u_2,\ldots,u_{s+1})\cr  
& + \sum_{1\leq k\leq s} (-1)^k C(u_1,\ldots,u_k u_{k+1},\ldots,u_{s+1})\cr
& +(-1)^{s+1} C(u_1,\ldots,u_s) u_{s+1}.
\end{align*}
The $sth$ Hochschild cohomology space will be denoted by $H^s({\mathcal A},M)$ or 
$H{}_{diff}^s({\mathcal A},M)$ if we restrict
ourselves to differential cochains.}
 
\vskip0,5cm
Now, let ${\mathcal C}{}_n^s=C{}_{n,grad,nc}^s(S(\u{g}))$ be the space of homogeneous of degree $-n$, correct of degree $-n$,
vanishing on constants, $s$-linear operators on $S(\u{g})$, differential or not. We denote by 
${\mathcal C}_{n,diff}^s=C{}_{n,grad,nc,diff}^s(S(\u{g}))$ the subspace of operators of ${\mathcal C}{}_n^s$ which are
differential. 
  $({\mathcal C}{}_n^\ast,\d)$ and $({\mathcal C}{}_{n,diff}^\ast,\d)$ are  subcomplexes
 of $(C^\ast(S(\u{g}),S(\u{g})),\d)$.
These subcomplexes give rise to 
well graded cohomology spaces  denoted by $H{}_{n,grad,nc}^{\ast}(S(\u{g}))$
and $H{}_{n,grad,nc,diff}^{\ast}(S(\u{g}))$. 
It is useful to know when these spaces vanish. In particular, we have the following result.
\par\noindent
{\bf Proposition 7.}\par\noindent
{\sl
   $\displaystyle H{}_{n,grad,nc}^3(S(\u{g}))=H{}_{n,grad,nc,diff}^3(S(\u{g}))=\{0\}\quad \forall n
\geq 4$ .}

{\sl Proof:}
1) Let $E$ be an element of ${\mathcal C}{}_{n,grad,nc}^3(S(\u{g}))$ such that $\d E=0$.
Clearly, $E(u,v,w)$ can be decomposed into a sum of two cocycles $ E_1 + E_2$ with $E_1$
symmetric in $u,w$ and  $E_2$ skew-symmetric in $u,w.$
In [14] p.148, G. Pinczon shows that if $N$ denotes the algebra of smooth functions over $\u{g}^\ast$,
then
$$\displaystyle H^3(S(\u{g}), N)= H{}_{diff}^3(S(\u{g}),N)= H{}_{diff}^3(N,N).$$
 It is well known $[16]$ that $H_{diff}^{\ast}(N,N)$ is isomorphic to the space of skew multivectors fields over
$\u{g}^\ast$. Thus, there exist two Hochschild cochains $C_1$ and $C_2$ in $C^2(N,N)$ such that

$\bullet E_1=\d C_1$ with $C_1(u,v)$ skew-symmetric in $u,v$

$\bullet E_2=\d C_2 + A$ with $C_2(u,v)$ symmetric in $u,v$ and
where $A(u,v,w)$ is 
$$A(u,v,w)=\sum a_{ijk} \p_i(u) \p_j(v) \p_k(w)$$ with completely skew-symmetric coefficients.\par
\noindent Since
$$\displaystyle \oint_{(u,v,w)} E_2(u,v,w)
:= E_2(u,v,w)+E_2(v,w,u)+E_2(w,u,v)=3\,A(u,v,w),$$
\noindent $A$ is necessarily homogeneous of degree $-n$ and the coefficients $a_{ijk}$ are polynomials of degree 
$3-n$. Thus, for $n\geq 4$, $n-3<0$ and $A\equiv 0.$ 
Moreover, since $E$ vanishes on constants, we can suppose it is the same for $C_1$ and $C_2$ just by replacing
$C_i$ by $C_i -\d T_i$ where $T_i$ is defined by $T_i(u)=C_i(u,1)$  ($i=1,\,2$).  
\vskip0,5cm\noindent
2) Let us now prove that $C_1$ can be chosen in ${\mathcal C}{}_{n,grad,nc}^2(S(\u{g}))$:
$E_1$ and $C_1$ can be decomposed into an infinite sum 
$\displaystyle \sum_{ N\geq 0} E_{1,N}\quad(\sum_{N\geq 0} C_{1,N}$ respectively) of differential operators in the 
variables $(x_i)$ of the form
$$ E_{1,N}(u,v,w)=\displaystyle\sum_{a+b+c=N} E_{k_1...k_a,l_1...l_b,m_1,...,m_c}
\p_{k_1,...,k_a}(u)
\p_{l_1,...,l_b}(v)\p_{m_1...m_c}(w).$$ Respectively,
$$\displaystyle
 C_{1,N}(u,v) =\sum_{a+b=N,a\geq b}C_{k_1...k_a,l_1...l_b}
(\p_{k_1...k_a}(u)\p_{l_1...l_b}(v)-
\p_{k_1...k_a}(v)\p_{l_1...l_b}(u))\,(*)$$
where the coefficients $C_{k_1...k_a,l_1...l_b}$ are supposed to be symmetric in the indexes $k_i$ and in
the  indexes $l_j$, and such that $C_{k_1...k_a,l_1...l_b}=-C_{l_1...l_b,k_1...k_a}$ (if
$a=b$).\par\noindent Since $E_1$ vanishes on constants, $E_{1,N}=0$ if $N<3$. 
Thus, $$E_1=\displaystyle \sum_{N\geq 3} E_{1,N}
:=\sum_{N\geq 3} (E_1)_N=\sum_{N\geq 3} (\d C_1)_N=\sum_{N\geq 3} \d(C_{1,N}).$$
The last equality directly comes from the definition of the Hochschild coboundary. In the following, we shall
assume that $C_1=\displaystyle \sum_{N\geq 3} C_{1,N}$ because $C_{1,1}$ and $C_{1,2}$ are not involved in the expression
of $E_1=\d C_1$.\par
Then, we want to prove that every $C_{1,N}\,(N\geq 3)$ sends $S(\u{g})\times S(\u{g})$ into $S(\u{g})$ and is 
homogeneous and correct of degree $-n$, or equivalently, to show that every 
$C_{k_1...k_a,l_1...l_b}$ is an
element of $S(\u{g})$ homogeneous of degree $a+b-n=N-n$ and that 
$$|| C_{k_1...k_a,l_1...l_b}||\leq ||X_{k_1}||+\ldots+|| X_{k_a}||+||X_{l_1}||+\ldots+||X_{l_b}||-n.$$
To this end, we  use a technique which can be found in [11] p.238-242 or in 
Gutt's thesis [9].\par\noindent
By $(*)$, $C_{1,N}$ is a finite sum of terms of type $(a,b)\,(a\geq b)$. Let $(r,s)$ be the highest of the types $(a,b)$
($(a,b) >(a',b')$ if $\,a >a'$ or if $\,a=a'$ and $\,b>b'$).
Let also $C_{k_1...k_r,l_1...l_s}$ be the topmost coefficient with respect to lexicographical order in
the indexes.  We shall call principal part of $C_{1,N}$ the unique term $P$ of type $(r,s)$ defined by 
$$\displaystyle P(u,v) =C_{k_1...k_r,l_1...l_s}(\p_{k_1...k_r}(u)\p_{l_1...l_s}(v)-
\p_{k_1...k_r}(v)\p_{l_1...l_s}(u)).$$ The principal part $P$ of $C_{1,N}$ becomes in $\d (C_{1,N})$
$$-C_{k_1...k_r,l_1...l_s}(\p_{k_1...k_r}(uv)\p_{l_1...l_s}(w) -\p_{k_1...k_r}(w)
\p_{l_1...l_s}(uv))$$
$$+C_{k_1...k_r,l_1...l_s}(\p_{k_1...k_r}(u)\p_{l_1...l_s}(vw)-\p_{k_1...k_r}(vw)\p_{l_1...l_s}(u))$$
up to terms without any derivatives of $u$ or of $w.$ \par\noindent
By construction, there are only three cases to consider:\par 
$\circ$ if $r\geq s, r\geq 2, s\geq 2:$
There is only one term of type $(r,s-1,1)$ in $\d (C_{1,N})$ which corresponds to the principal part and which can
be written, up to some constant coefficient,
$$\displaystyle C_{k_1...k_r,l_1...l_s}\p_{k_1...k_r}(u)\p_{l_1...l{s-1}}(v)\p_{l_s}(w).$$

$\circ$ if $r\geq 3, s=1:$
The only term of type $(r-1,1,1)$ in $\d (C_{1,N})$ can be written, up to an eventual constant coefficient,
$$\displaystyle C_{k_1... k_r,l_1}\p_{k_1... k{r-1}}(u)\p_{ k_r}(v)\p_{l_1}\,(w).$$
In these two cases, the coefficients $C_{k_1... k_r,l_1... l_s}$ are convenient.

$\circ$ if $r=2, s=1:$
We get the following terms of type $(1,1,1)$ in $\d (C_{1,N})$, up to some constant coefficient,
$$\displaystyle (C_{ij,k}+ C_{jk,i})\p_i(u)\p_j(v)\p_k(w).$$ 
Thus, $(C_{ij,k}+C_{jk,i})$ is polynomial of degree $3-n$. For $n\geq 4$, $n-3<0$ thus
$C_{ij,k}+C_{jk,i}=0.$
By cyclic summation, we find $C_{ij,k}=0.$ In other words, there are no terms of type $(2,1)$ in $C_{1,N}$.
Then, as we see, the principal part $P$ of $C_{1,N}$ is homogeneous and correct of degree $-n$.
We can now repeat the proof for the principal part of  $(C_{1,N}-P)$... 
A step-by-step application of the same arguments finally shows that all the $C_{1,N}$ (thus also $C_1=\sum C_{1,N}$) 
belong to ${\mathcal C}{}_n^2.$
\vskip0,5cm \noindent
3) We apply the same method for $C_2$: we start by decomposing $E_2$   
$$\displaystyle E_2=\sum_{N\geq 3} E_{2,N}=\sum_{N\geq 3} (\d C_2)_{N}=\sum_{N\geq 3} \d(C_{2,N}).$$
We can suppose that $C_2=\displaystyle \sum_{N\geq 3} C_{2,N}.$\par\noindent
Now, the principal part of $C_{2,N}$ of type $(r,s)$ can be written as follows
$$C_{k_1,..., k_r,l_1,...,l_s}(\p_{k_1,...,k_r}(u)\p_{l_1,...,l_s}(v)+\p_{k_1,...,k_r}(v)
\p_{l_1,...,l_s}(u)).$$

$\circ$ The cases $r\geq s\geq 2$ and $r\geq 3, s=1$ are the same as above.

$\circ$ If $r=2, s=1:$
The terms of type $(1,1,1)$ in $\d(C_{2,N})$, up to some constant coefficients, are
$$(C_{ij,k}-C_{jk,i}) \p_i(u)\p_j(v) \p_k(w).$$
Since $n\geq 4$, $(C_{ij,k}-C_{jk,i})$ is polynomial of degree $n-3<0$ and $C_{ij,k}=C_{jk,i}$. Therefore, the $C_{ij,k}$
do not appear in $\d (C_{2,N})$
and we can  remove every terms of type $(2,1)$ from the expression of $C_{2,N}$. As before, we succeed in proving that
all the  $C_{2,N}$ with $N\geq 3$ (thus also $C_2=\displaystyle \sum_{N\geq 3} C_{2,N}$) are  elements 
of ${\mathcal C}_{n,grad,nc}^2(S(\u{g}))$.\par\noindent
This ends the proof for non-differential cochains. Obviously, the same can be done for differential cochains. 
\par\medskip\noindent
{\bf Remark 8.}\par\noindent
{\sl If $E_1$ is a symmetric $3$-cocycle in ${\mathcal C}{}_3^3$, 1) and 2) are
still valid. There exists $C_1$ in
${\mathcal C}{}_3^2$ such that $E_1=\d C_1$. The only difference is that the terms of type $(2,1)$ of $C_1$ have now constant
coefficients.
}

\vskip1cm

\section{Tangential and well graded deformation of $S(\u{g})$}\label{sec3}

\smallskip
\noindent
{\bf Definition 9.}\par\noindent
{\sl A graded star product of $S(\u{g})$ is a bilinear map from $S(\u{g})\times S(\u{g})$ to
$S(\u{g})[[\nu]]$ defined by 
$$(u,v)\rightarrow \displaystyle u\ast v=uv+\{u,v\}\nu +\sum_{n\geq2} C_{n}(u,v){\nu}^n $$
where the cochains $C_n$ are operators on $S(\u{g})$ with values in $S(\u{g})$ satisfying the following properties
\par\noindent
\begin{align*}
&(i)\,C_n(u,v)=(-1)^n C_n(v,u),\forall u,v\in S(\u{g});\cr
&(ii)\,C_n(1,v)=0,\forall v\in S(\u{g});\cr
&(iii)\,C_n\,\,is\,\,homogeneous\,\,of\,\,degree\,-n; \cr
& (iv) \displaystyle \sum_{r+s=k}C_r(C_s(u,v),w)=\sum_{r+s=k}C_r(u,C_s(v,w)), 
k\geq 0,\forall u,v,w\in S(\u{g}).
\end{align*}}

\noindent
{\sl $\ast$ is said to be a well graded star product of $S(\u{g})$ if the cochains $C_n$ are both
homogeneous of degree
$-n$ and correct of degree $-n$.}

\noindent
{\sl Moreover, $\ast$ is said to be a tangential star product of $S(\u{g})$ if the $C_n$ are tangential
operators on
 $S(\u{g})$.}
\par\medskip 

\noindent
{\bf Theorem 10.}\par\noindent
{\sl Suppose we know a tangential operator $C_2$ on $S(\u{g})$ homogeneous and correct of degree $-2$,
such that 
$\displaystyle u\ast v=uv + \{u,v\}\nu + C_2(u,v)\,\nu^2$
is associative up to order $3$ in $\nu$. Then $C_2$ is the second order term of a		tangential, well graded 
star product of $S(\u{g})$. Moreover, if $C_2$ is differential, $C_2$ is the second order term of a
differential,  tangential and well graded star product of $S(\u{g})$.}\par
{\sl Proof:} Let us assume that we have found tangential, homogeneous and correct
$C_2,\ldots,C_{n-1}$ $(n\geq 3)$ such that
$u\ast v=uv + \{u,v\}\nu + \displaystyle \sum_{2\leq k\leq n-1} C_k(u,v)\nu^k$ is associative up to order $n$ in $\nu$.
Consider then the Hochschild cocycle $E_n$ defined by
$$E_n(u,v,w)=\displaystyle \sum_{r\geq 1,s\geq 1,r+s=n} C_r(C_s(u,v),w)-C_r(u,C_s(v,w)).$$
Clearly, $E_n$ is homogeneous and correct of degree $-n$. Thanks to Proposition 7 and Remark 8, we can
already say that there exists an operator $C_n$ on $S(\u{g})$ such that 
$E_n=\d C_n$, $C_n(u,v)=(-1)^n C_n(v,u)$, $C_n(1,v)=0$ for all $u,v$ in $S(\u{g})$ and so that $C_n$ is
both homogeneous of degree $-n$ and correct of degree $-n$. It remains to show that $C_n$ is tangential. \par\noindent
First, by transposing the equality $E_n=\d C_n$ in coordinates $(p,q,\l)$, one obtains:
 $\widehat{E_n}=\d ({\widetilde{C_n}}_{(1)})$
where $\widehat{\,}$ and $\widetilde{\,}$ have the same meaning as in Section $1$. To make the writing
simpler, we forget the $n$ and we note $E=E_n$ and $C=C_n$. As in Proposition 7, we decompose $\hat{E}$ and
${\widetilde{C}}_{(1)}$ in an  infinite sum 
$\displaystyle \sum_{K\geq 0} {\hat{E}}_K\quad(\sum_{K\geq 0} {\widetilde{C}}_K$ respectively) of operators in the variables 
$(p,q,\l)$ of the form
$${\hat{E}}_K (u,v,w)=\displaystyle \sum_{a+b+c=K} E_{k_1... k_a,l_1...l_b,m_1...m_c}
\p_{k_1...k_a}(u) \p_{l_1...l_b}(v) \p_{m_1...m_c}(w).$$
Respectively,
$${\widetilde{C}}_K(u,v)=\displaystyle \sum_{a+b=K, a\geq b} C_{k_1...k_a,l_1...l_b}
(\p_{k_1...k_a}u\p_{l_1...l_b}v +(-1)^n \p_{k_1...k_a}v\p_{l_1...l_b}u).$$
 $E$ vanishes on constants, $\hat{E}_K=0$ for $K<3$. Thus, 
$$\displaystyle \sum_{K\geq 3} {\hat{E}}_K=\sum_{K\geq 3} \d({\widetilde{C}}_K).$$\noindent
{\bf First case:}  $n$ is odd $(n\geq 3)$\par\noindent
${\widetilde{C}}_{(1)}=\displaystyle\sum_{K\geq 3} \widetilde{C}_K + {\widetilde{C}}_2.$ But, since $C$ is correct of degree $-n$  
and since $n\geq 3$, ${\widetilde{C}}_2=0.$
Let us now prove that all the ${\widetilde{C}}_K\,(K\geq 3)$ do not involve derivatives with respect to $(\l_k)$.
To this end, we proceed
as in Proposition 7. We consider first the principal part $P$ of ${\widetilde{C}}_K$ of type $(r,s)$ in the
variables
$(p,q,\l)$. Three cases have to be considered. The cases $r\geq s\geq 2$ and $r\geq 3,s=1$ are directly solved. Now, if
$r=2$ and $s=1$, the terms of type $(1,1,1)$ in $\d \widetilde{C}_K$ up to some constant coefficients are
$$(C_{ij,k} - C_{jk,i})\p_i(u)\p_j(v)\p_k(w).$$
Suppose that some derivative with respect to $(\l_k)$ appears, then $(C_{ij,k}-C_{jk,i})$ should be zero. And, by
cyclic summation, we find $C_{ij,k}=0.$
Therefore, we conclude that there are no derivatives with respect to $(\l_k)$ in the principal part.
Then, we repeat the proof step by step and finally get that all the ${\widetilde{C}}_K$ (thus also
${\widetilde{C}}_{(1)}=\displaystyle \sum_{K\geq 3}{\widetilde{C}}_K$) do not involve derivatives with respect to the
variables $(\l_k)$.
Thus, $C=C_n$ is tangential if $n$ is odd.\par\noindent
{\bf Second case:} $n$ is even $(n\geq 4)$\par\noindent
${\widetilde{C}}_{(1)}=\displaystyle\sum_{K\geq 4} \widetilde{C}_K + {\widetilde{C}}_3 + {\widetilde{C}}_2.$ But, since $C$ is correct
of degree $-n$ and $n\geq 4$, ${\widetilde{C}}_2={\widetilde{C}}_3=0.$
As before, we use principal parts to show that all the $\widetilde{C}_K$ for $K\geq 4$ do not include derivatives
with respect to $(\l_k)$. In other words, $C=C_n$ is also tangential if $n$ is even. This ends the proof.

\medskip\noindent
{\bf Remark 11.}\par\noindent It is possible to define cohomology spaces related to tangential and well
graded deformations of
$S(\u{g})$. For the moment, we denote by ${\mathcal C}{}_{n,tang,grad}^s$ the space
of $s$-linear operators on $S(\u{g})$,  which are tangential, homogeneous of degree $-n$ and correct of degree $-n$.
Endowed with the Hochschild coboundary, ${\mathcal C}{}_{n,tang,grad}^{\ast}$
becomes a complex. Let $H{}_{n,tang,grad}^{\ast}(S(\u{g}))$ be the corresponding cohomology.
Then, we see that Proposition 7 together with Theorem 10 contain the computation of 
$H{}_{n,tang,grad}^3(S(\u{g}))$. Similarly, one can prove the vanishing of the second spaces of this
cohomology, $H{}_{n,tang,grad}^2(S(\u{g}))$, for
$n\geq 2$. More exactly, the following facts hold\par
$(i)$  if $C(u,v)$ is a cocycle, skew-symmetric in $u,v$, tangential, homogeneous of degree $-(2k-1)$ and correct
of degree $-(2k-1)$, 
$k\geq2$, then $C\equiv0;$

$(ii)$ if $C(u,v)$ is a cocycle, symmetric in $u,v$, tangential, homogeneous of degree $-2k$ and correct 
of degree $-2k$, $k\geq1$,
then we can suppose that $C=\d R$ where $R$ is tangential, homogeneous of degree $-2k$ and correct of degree $-2k$.

\medskip\noindent
{\bf Theorem 12.}\par\noindent
{\sl Two tangential and well graded star products of $S(\u{g})$ are always tangentially equivalent.
One can find an equivalence operator of the form
$$\displaystyle T=Id + \sum_{k\geq1} T_{2k}\nu^{2k}$$
where all the $T_{2k}$ are tangential operators from $S(\u{g})\times S(\u{g})$ to $S(\u{g})$, homogeneous of degree $-2k$
and correct of degree $-2k$.
(The homogeneity property also implies that for all $k$, each term of $T_{2k}$ is of order $\geq 2k$.)}

{\sl Proof:}
The result is a straightforward consequence of the previous remark.
Indeed, let $\ast$, $\ast'$ be two tangential and well graded star products of $S(\u{g})$. Let $k$ be
$\geq 1$. 
 Assume that we found $T_0,\ldots,T_{2k-2}$, such that $T_0=Id$, that every $T_{2j}\,(j\geq1)$ is tangential, 
 homogeneous of degree $-2j$, correct of degree $-2j$ and that the star product $\ast''$ defined by
$$\displaystyle u\ast''v=H^{-1}(H(u)\ast' H(v))$$
where 
$ H= Id + \ldots + T_{2k-2}{\nu}^{2k-2}$, satisfies 
$C_j''(u,v)=C_j(u,v)$ for all $j\leq{2k-2}.$
The associativity condition leads to
$$ \d(C_{2k-1}''-C_{2k-1})=0.$$
Afterwards, either $k=1$ and $C_1''(u,v)=C_1(u,v)=\{u,v\}$ or $k\geq2$ and $C_{2k-1}''=C_{2k-1}$ (Remark 11, (i)) .
Then, we obtain:
$$ \d(C_{2k}''-C_{2k})=0.$$
Thus, there exists $T_{2k}$ as announced (Remark 11, (ii)) so that $C_{2k}''=C_{2k}+\d{T}_{2k}.$
A simple induction enables us to construct the equivalence operator $T$ and thereby ends the proof.
\vskip1cm

\section{Applications and examples}\label{sec4}

\smallskip
Let us first recall the construction of the Gutt star product $\ast_{G}$ defined on the symmetric algebra
$S(\u{g})$ of any Lie algebra $\u{g}$. Let $U(\u{g})$ be the universal enveloping algebra of $\u{g}$
and $\s:S(\u{g})\rightarrow U(\u{g})\,$  be the symmetrization map. Denote by
$[\,u\,]_k$ the $kth$ component of an element $u$ of $U(\u{g})$ relative to the canonical decomposition
$\displaystyle U(\u{g})=\oplus \s(S^k(\u{g}))$. If $P,\,Q$ are homogeneous
polynomials of degree $r,\,s$ respectively, then 
$$\displaystyle P{\ast}_G Q=\sum_{n\geq0} C_{n,G}(P,Q) {\nu}^n=\sum_{n\geq0}{\s}^{-1}([\s(P).\s(Q)]_{(r+s-n)})\,(2\nu)^n.$$
Using linearity to extend the above expression to all polynomials, we get
the Gutt star product. In the literature, the same star product is sometimes called the star product
coming from the enveloping algebra via Poincar\'e-Birkhoff-Witt. One checks that $\ast_{G}$ is
differential and graded. \par
 Now, we are interested in the example of  $\u{g}_{54}$. This
nilpotent Lie algebra is defined by the following brackets
$$[X_5,X_4]=X_3, [X_5,X_3]=X_2, [X_4,X_3]=X_1.$$ 
The quotient field of $I(\u{g}_{54})$ is generated by two central elements, namely $X_1$ and $X_2$, and by
$\Delta\,=\,\frac{X_3^2}{2}\,+\,X_1X_5\,-\,X_2X_4.$

\par
A simple calculation shows that, up to a normalization, the second order term $C_{2,G}$ satisfies
$$\displaystyle
C_{2,G}(\Delta,.)=\frac{x_1^2}{6}\p_{44}+\frac{x_1x_2}{3}\p_{45}+\frac{x_2^2}{6}\p_{55}.$$
 It
is thus clear that $\ast_G$ is not tangential. If it was, $C_{2,G}(\Delta,.)$ would be reduced to
zero.\par
 A natural idea to know whether a tangential and well graded star product of
$S(\u{g}_{54})$ exists or not, is to try to correct $C_{2,G}$ by means of an operator $T$ on
$S(\u{g}_{54})$ such that
$\displaystyle C_{2,G}(\Delta,v)+\d T(\Delta,v)=0$ for all $v$ in $S(\u{g}_{54}).$
A possible $T$ is
$$
\displaystyle T\, =\sum_{n\geq 4} (-1)^{n} \frac{2^{n-3}}{6 (n-2)!}\, x{}_3^{n-4}
 \Big( x_2^2 \p_3^{n-2}\p_{55}
  + x_1^2 \p_3^{n-2} \p_{44} + 2 x_1 x_2 \p_3^{n-2} \p_{45}\Big).$$\noindent
If we note $\,{\s}_3=\displaystyle \sum_{n\geq0}\frac{(-2x_3)^n}{n!}{\p_3}^n$, $T$ can be written in
the form  \par\noindent
$\displaystyle T\,=\,A_{55} \p_{55}+A_{45} \p_{45}+A_{44} \p_{44}+A_{355} \p_{355}+A_{345} \p_{345}+A_{344} \p_{344},$
where

\begin{align*}
A_{55}\,=\,&\frac{{x_2}^2}{12 {x_3}^2}({\s}_3-Id) \cr
A_{45}\,=\,&\frac{x_1x_2}{{x_3}^2}({\s}_3-Id) \cr
A_{44}\,=\,&\frac{{x_1}^2}{12{x_3}^2}({\s}_3-Id) \cr
A_{355}\,=\,&\frac{{x_2}^2}{6x_3}\cr
A_{345}\,=\,&\frac{{x_1x_2}}{3x_3}\cr
A_{344}\,=\,&\frac{{x_1}^2}{6x_3}.
\end{align*}
One immediately sees that $C_2=C_{2,G}+\d T$ is tangential and homogeneous of degree $-2$. 
Now, we need to prove that $C_2$ is also correct of degree $-2$.
Let us first introduce the canonical variables 
$$(p=x_4,q=\frac{x_3}{x_1}, \l_1=x_1, \l_2=x_2, \l_3=(\frac{{x_3}^2}{2}+x_1x_5-x_2x_4)).$$
Then, changing variables ($(x_i)\rightarrow (p,q,\l)$) and using the notation of Section 1, we obtain
\begin{align*}
 \widetilde{C_{2,G}}(u,v)=&{\p}_{p p}(u){\p}_{q q}(v)-2{\p}_{p q}
(u){\p}_{p q}(v)+{\p}_{q q}(u){\p}_{p p}(v)\cr
&+\frac{1}{3}\l{}_1^2(\p_{\l_3 p}(u)\p_p(v)+\p_{\l_3 p}(v)\p_p(u))\cr
&-\frac{1}{3}\l{}_1^2(\p_{pp}(u)\p_{\l_3}(v)+\p_{pp}(v)\p_{\l_3}(u))\cr
\widetilde{T}_{(1)}(u)=&\displaystyle\sum_{n\geq 4}(-1)^n
\frac{2^{n-3}}{6(n-2)!} q^{n-4}(q \,{\l}_1^2 {\p}_{\l_3} + {\p}_q)^{n-2} \p_{pp}(u).
\end{align*}
 Recall now that $\widetilde{C_2}_{(1)}=\widetilde{C_{2,G}}+\d \widetilde{T}_{(1)}$ coincides with $C_2$ on 
$S(\u{g})$. Thus, the above expressions  mean that $C_2$ is correct of degree $-2$.
 By Theorem $10$ (see Section $3$), this is sufficient to show the existence of an 
algebraic, tangential and well graded star product on $\u{g}_{54}^\ast$. That is the best we can do,
because there is no deformation on $\u{g}_{54}^\ast$ which is both tangential and differential [1, 5]. 
Remark also that, since the only polynomials on $\u{g}{}_{54}^\ast$, whose restriction to an orbit of maximal dimension 
is zero, are invariant (see [13] p.23), our deformation is tangential to all the regular orbits ({\sl i.e.} orbits of
maximal dimension).

\par\medskip

Nevertheless, as differential operators are more convenient to handle than algebraic maps, let us mention the possibility
of constructing a differential and tangential deformation on the subset $\Omega$ of all the orbits of maximal dimension.
Note that $\Omega$ is a regular Poisson manifold and that    
$$\Omega=\{ \xi=({\xi}_1,...,{\xi}_5) 
\in \u{g}{}_{54}^{\ast}\quad\hbox{such that}\quad\e{}_1^2+ \e{}_2^2
+\e{}_3^2\not=0\}.$$  We found an explicit expression of an operator $C_2^{\k}$ with homogeneous
coefficients in
$C^{\infty}(\Omega)$, which is both tangential and differential. Here it is
$$ C_2^{\k}=C_{2,G} +\d T^{\prime},$$
where
\begin{align*}
T'\,=\,& A_{453}\p_{453}+A_{355}\p_{355}+A_{455}\p_{455}+A_{344}\p_{344}+A_{445}\p_{445}\cr
    +\,&A_{555}\p_{555}+A_{444}\p_{444}\cr
r\,=\,&x_1^2+x_2^2+x_3^2 \cr
A_{453}\,=\,&\frac{x_1x_2x_3}{3r}\cr
A_{355}\,=\,&\frac{x_3 x_2^2}{6r}\cr
A_{455}\,=\,&\frac{-x_2^3+2x_1^2x_2}{6r}\cr
A_{344}\,=\,&\frac{x_1^2x_3}{6r}\cr
A_{445}\,=\,&\frac{x_1^3-2x_1x_2^2}{6r}\cr
A_{555}\,=\,&\frac{x_1x_2^2}{6r} \cr
A_{444}\,=\,&\frac{-x_1^2x_2}{6r}.\end{align*}
\smallskip
Further examples are given by Pedersen in [13]. Let us say a few words about 
$\u{g}_{6\,12}$ ([13] p.87) and 
$\u{g}_{6\,14}$ ([13] p.99).\par
$\bullet$ The Lie algebra structure of $\u{g}_{6\,12}$ is defined by the non vanishing brackets
$$[X_6,X_5]=X_4, [X_6,X_4]=X_3, [X_6,X_3]=X_2, [X_5,X_2]=-X_1, [X_4,X_3]=X_1.$$ 
The quotient field of $I(\u{g}_{6\,12})$ is generated by 
$X_1$ and $\frac{X_3^2}{2}-X_2X_4+X_1X_6$.\par
$\bullet$ The Lie algebra structure of $\u{g}_{6\,14}$ is defined by the following brackets
$$[X_6,X_5]=X_4,[X_6,X_4]=X_3,[X_6,X_3]=X_2,
[X_5,X_4]=X_2,[X_5,X_2]=-X_1,$$\par\noindent$$[X_4,X_3]=X_1.$$
Moreover, the quotient field of $I(\u{g}_{6\,14})$ is generated by
$X_1$ and $\frac{X_2^3}{3}-\frac{X_1X_3^2}{2}+X_1X_2X_4-X_1^2X_6$.\par\noindent
For these two examples, we may explicitly define an algebraic, tangential $C_2$ on the symmetric algebra, and 
a differential, tangential $C_2'$ on $C^{\infty}(\Omega)$, $\Omega$ denoting the open set of regular orbits, with similar 
argument as for $\u{g}_{54}$.
\vskip1cm
\centerline{\bf Acknowledgements}
\smallskip
I am very deeply grateful to D. Arnal for his constant help and for so much fruitful advice.
I also want to thank M. Masmoudi for some interesting remarks.
The research has been supported by an ARC of the Communaut\'e Fran\c{c}aise de Belgique.
\vskip1cm
\centerline{\bf References}
\smallskip

\noindent
1. Arnal, D., Cahen, M. and Gutt, S.: Deformations on coadjoint orbits, {\sl J. Geom. Phys.}
{\bf Vol 3 n.3}  (1986).

\noindent
2. Asin Lares, S.: On tangential properties of the Gutt star-product, {\sl J. Geom. Phys.} {\bf 24}
(1998), 164-172.

\noindent
3. Bayen, F., Flato, M., Fronsdal, C., Lichnerowicz, A. and Sternheimer, D.: Deformation theory and quantization, 
{\sl Ann. Phys.} {\bf III} (1978), 61-151.

\noindent
4. Bonnet, P.: Param\'etrisation du dual d'une alg\`ebre de Lie nilpotente, {\sl Ann. Inst. Fourier, Grenoble}
{\bf 38}, 3 (1988), 169-197.

\noindent
5. Cahen, M., Gutt, S. and Rawnsley, J.: On tangential star-products for the coadjoint Poisson structure,
{\sl Comm. Math. Phys.} {\bf 180} (1996), 99-108.

\noindent 
6. Cahen, M. and Gutt, S.: An algebraic construction of $\ast$-product on the regular orbits of semi-simple Lie groups.
Gravitation and geometry, 71-82, Monographs Textbooks Phys. Sci., {\bf 4}, Bibliopolis, Naples (1987). MR 89e:22023.

\noindent
7. De Wilde and Lecomte, P.B.A.: Existence of star-products and formal deformations in Poisson Lie algebra of 
arbitrary symplectic manifold, {\sl Lett. Math. Phys.} {\bf 7} (1983), 487-496.

\noindent
8. Fedosov, B.V.: A simple geometrical construction of deformation quantization, {\sl J. Diff. Geom.}
{\bf40} (1994), 213-238.

\noindent
9. Gutt, S.: D\'eformations formelles de l'alg\`ebre des fonctions diff\'erentiables sur une vari\'et\'e
symplectique, Th\`ese Bruxelles (1980).

\noindent 
10. Kontsevich, M.: Deformation quantization of Poisson manifolds, I Preprint q-alg/9709040 (1997).

\noindent
11. Lichnerowicz, A.: Vari\'et\'es de Poisson et feuilletages, {\sl Ann. Fac. Toulouse} (1982), 195-262.

\noindent
12. Masmoudi, M.: Star-produits sur les vari\'et\'es de Poisson, Th\`ese, {\sl Universit\'e de Metz} (1992).

\noindent
13. Pedersen, N.V.: Geometric quantization and nilpotent Lie groups. A collection of examples, {\sl University of
Copenhagen Denmark} (1988), 1-180.

\noindent
14. Pinczon, G.: On the equivalence between continuous and differential deformation theories, {\sl Lett. Math. Phys.} {\bf39}
(1997), 143-156.

\noindent
15. Vergne, M.: La structure de Poisson sur l'alg\`ebre sym\'etrique d'une alg\`ebre de Lie nilpotente,
{\sl Bull. Soc. Math. Fr.} {\bf 100} (1972), 301-335.

\noindent
16. Vey, J.: D\'eformation du crochet de Poisson sur une vari\'et\'e symplectique, {\sl Comment. Math. Helv.} {\bf 50} (1975), 
421-454. 

\end{document}